\documentclass[10pt,twocolumn,numbers]{article}

\usepackage{graphicx}

\usepackage{amssymb,amsthm,amsmath,dsfont,bbm,mathtools,stmaryrd}

\usepackage{lineno}

\usepackage[colorlinks=true, citecolor=blue]{hyperref}
\usepackage[lined,ruled,commentsnumbered]{algorithm2e}
\usepackage{algcompatible}

\usepackage{wrapfig}
\usepackage[numbers,sort&compress]{natbib}
\usepackage{comment}
\usepackage{longtable}
\usepackage{multirow}
\usepackage{bbding}
  
\usepackage[bottom]{footmisc}
\usepackage{booktabs}
\usepackage{pifont}
\usepackage{float}
\usepackage{wrapfig}
\usepackage{framed,xcolor}
\usepackage{newfloat}
\usepackage[xcolor,framemethod=TikZ]{mdframed}
\usepackage{amsthm}
\usepackage{thmtools}
\usepackage{authblk}
\patchcmd{\thebibliography}{\section*{\refname}}{}{}{}
\usepackage{tikz}
\usetikzlibrary{arrows,shapes,chains,calc,patterns,decorations.pathreplacing}
\usepackage{stmaryrd}
\usepackage{fancyhdr}
\usepackage{subfigure}
\usepackage{comment}
\def\BibTeX{{\rm B\kern-.05em{\sc i\kern-.025em b}\kern-.08em
		T\kern-.1667em\lower.7ex\hbox{E}\kern-.125emX}}

\renewcommand{\headrulewidth}{2pt}

\newlength\FHoffset

\fancyheadoffset{\FHoffset}

\newlength\FHleft
\newlength\FHright

\newbox\FHline
\setbox\FHline=\hbox{\hsize=\paperwidth%
	\hspace*{\FHleft}%
	\rule{\dimexpr\headwidth-\FHleft-\FHright\relax}{\headrulewidth}\hspace*{\FHright}%
}

\newtheoremstyle{theoremdd}
{\topsep}
{\topsep}
{\itshape}
{0pt}
{\fontfamily{cmss}\selectfont\bfseries}
{.}
{ }
{\thmname{#1}\thmnumber{ #2}\thmnote{ (#3)}}

\theoremstyle{theoremdd}

\mdfdefinestyle{myStyle}{roundcorner=5pt,backgroundcolor=brown!20,linecolor=red!40!black,linewidth=1.5pt}

\usepackage{titlesec}
\titleformat*{\section}{\fontfamily{cmss}\selectfont\large\bfseries\color{red!40!black}}
\titleformat*{\subsection}{\fontfamily{cmss}\selectfont\normalsize\bfseries\color{red!40!black}}
\titleformat*{\subsubsection}{\fontfamily{cmss}\selectfont\normalsize\color{red!40!black}}
\usepackage[left=0.69 in, right=0.69 in, top=1.1 in, bottom=1.3 in]{geometry}
\graphicspath{{./Figures/}}

\newcommand\blfootnote[1]{%
	\begingroup
	\renewcommand\thefootnote{}\footnote{#1}%
	\addtocounter{footnote}{-1}%
	\endgroup
}
\renewcommand\abstractname{\fontfamily{cmss}\selectfont\normalsize\bfseries\color{red!40!black}\textbf{Abstract}}

\renewenvironment{abstract}{%
	\centering\small
	\list{}{\leftmargin1.5cm \rightmargin\leftmargin}
	\item\relax
	
	\begin{mdframed}[]
		\item[\hskip\labelsep\scshape\abstractname.]%
	}{%
	\end{mdframed}
	\endlist \par\bigskip
}
\makeatletter
\patchcmd{\@maketitle}{\LARGE \@title}{\fontfamily{cmss}\selectfont\LARGE\color{red!40!black}\@title}{}{}
\makeatother

\begin{document}

		
		
		\title{Generalized Adaptive Smoothing Using Matrix Completion for Traffic State Estimation}
		
			

\author[1]{Chuhan Yang$^{\star}$}
\author[1]{Bilal Thonnam Thodi}
\author[1,2]{Saif Eddin Jabari}

\affil[1]{New York University Tandon School of Engineering, Brooklyn NY, U.S.A.}
\affil[2]{New York University Abu Dhabi, Saadiyat Island, P.O. Box 129188, Abu Dhabi, U.A.E.}

\date{}


\twocolumn[
\begin{@twocolumnfalse}
	
\maketitle	

\begin{abstract}
	The Adaptive Smoothing Method (ASM) is a data-driven approach for traffic state estimation. It interpolates unobserved traffic quantities by smoothing measurements along spatio-temporal directions defined by characteristic traffic wave speeds. The standard ASM consists of a superposition of two a priori estimates weighted by a heuristic weight factor. In this paper, we propose a systematic procedure to calculate the optimal weight factors. We formulate the a priori weights calculation as a constrained matrix completion problem, and efficiently solve it using the Alternating Direction Method of Multipliers (ADMM) algorithm. Our framework allows one to further improve the conventional ASM, which is limited by utilizing only one pair of congested and free flow wave speeds, by considering multiple wave speeds. Our proposed algorithm does not require any field-dependent traffic parameters, thus bypassing frequent field calibrations as required by the conventional ASM. Experiments using NGSIM data show that the proposed ADMM-based estimation incurs lower error than the ASM estimation.
	
\end{abstract}
\bigskip
\end{@twocolumnfalse}
]

		
	
	
	

\section{INTRODUCTION}

Accurate\blfootnote{Corresponding author. Email: \url{cy1004@nyu.edu}.} knowledge of road traffic conditions plays a key role in real-time traffic management systems such as traffic lights, vehicle routing, and road performance evaluations \cite{li2021backpressure,li2019position,yang2016cav}. However, field traffic data obtained from the stationary detectors and floating cars, the two most popular traffic data sources, remain sparse in practice. One resorts to estimation techniques to infer missing traffic data from these sparse measurements.



A popular data-driven method for traffic state estimation is the Adaptive Smoothing Method (ASM), originally developed for stationary detector data \cite{treiber2002filter} and then later applied to floating cars data \cite{trieber2011filter}. ASM is an interpolation method based on the simplified kinematic wave theory of traffic flow \cite{newell1993simplified}. It assumes that perturbations in macroscopic traffic propagates forward (driving direction) in free-flow traffic and backward (in the opposite direction of traffic) in congested traffic. 
Accordingly, the ASM first builds two traffic estimates for congested and free flow traffic using an anisotropic low-pass filters. The final estimate is taken as a convex combination of these free-flow and congested traffic estimates using weights which depend on field traffic conditions.

Despite the simplicity and wide application of ASM, there are some limitations to be addressed. First, the ASM can only accommodate two traffic waves (one forward and one backward). In real-world traffic, a range of forward and backward traffic waves are observed. 
Second, the ASM weights used for combining free-flow and congested estimates is based on a heuristic formula that requires field calibration. Further, these weights are sensitive to near-capacity traffic conditions and produce inaccurate estimation results.

To address these shortcomings, we pose the weight calculation in the ASM as a constrained optimization problem, namely a kind of matrix completion problem. The proposed frameworks allows one to accommodate multiple traffic waves during estimation. The optimization problem can be efficiently solved using the alternating direction method of multipliers (ADMM). 
Our proposed method is calibrated with fewer hyperparameters and gives a better estimation error. Moreover, its performance can be further improved by including more a priori estimates, which were originally limited to a single congested and a single free flow estimate in traditional ASM implementations. 


The rest of the paper is organized as follows: Section \ref{sec2:relatedworks} briefly reviews the relevant literature. In Section \ref{sec3:methodology}, we present the ASM and our proposed modification using ADMM. The experiment details are provided in Section \ref{sec4:experiments}, followed by a brief discussion of results in Section \ref{sec5:results}. Finally we draw our conclusion and discuss on future research in Section \ref{sec6:conclusions}.


\section{Related Works}
\label{sec2:relatedworks}

Traffic state estimation techniques are broadly grouped into model-based,  data-driven learning and structured learning methods. Model-based methods combine estimations from a physical traffic flow model with field measurements using an exogenous filter \cite{papadopoulou2018microscopic,jabari2018stochastic,bekiaris2016tse_ieee,jabari2013gauss,hoogen2012lagrang,jabari2012stochastic}. The estimations from model-based methods are physically reasonable but are often limited by the capacity of traffic flow model and the filtering assumptions. Data-driven methods build parametric/non-parametric machine learning models from historic traffic data \cite{li2022ensemblemc,jabari2019learning,benkraouda2020traffic,yuhan2016dnn-speed,li2020short,xiao2018speed}. Exploiting non-linear regression functions such as deep neural networks or support vector machines, the data-driven techniques better capture the higher order traffic features and are more accurate than the model-based methods. But, lack of physical interpretability (black box nature) often limits these data-driven methods from practical applications.

The third group of estimation methods consist of structured-learning methods, where the data-driven methods are built that honor traffic physics constraints such as conservation laws and kinematic wave theory to improve the interpretability \cite{thodi2021anisocnn,thodi2021itsc,shi2021physics,huang2020physics,kaidi2019queueest}. The physical constraints are honored either during the model fitting stage \cite{shi2021physics,huang2020physics} or infused in the model architecture \cite{thodi2021anisocnn,jabarisparse2020,jabari2019learning,jabari2018stochastic}. These methods have shown robust estimation performance and requires limited data in the function fitting process.

We group the Adaptive Smoothing Method (ASM) \cite{treiber2002filter,trieber2011filter} in the category of structured-learning methods, since the two dimensional interpolation in ASM is not adhoc but takes into account the wave propagation characteristics in free-flow and congested traffic conditions. Different researchers tried to improve and modify the conventional ASM. For instance, \cite{schreiter2010two} improve the computational time of ASM using Fast Fourier Transform techniques. In another study, \cite{chen2019filter} propose to dynamically change the ASM kernel parameters (wave speeds and kernel smoothing widths) in a rolling-horizon framework. But the weights used for combining a priori traffic estimates of free-flow and shockwave wave speeds are unchanged and are based on a heuristic formula. The method proposed in our paper overcomes this limitation by systematically deriving the optimal weights for different wave speeds considered in ASM. The recent study from \cite{thodi2021anisocnn} uses the ASM anisotropic kernels in designing efficient convolutional neural networks for traffic speed estimation, where the wave speeds and their respective weights are learned from data.


\section{METHODOLOGY}
\label{sec3:methodology}

In this section, we first briefly describ the ASM, then we specify how to compute the weights using ADMM. 

\subsection{Notation}
We denote matrices using uppercase bold Roman ($\mathbf{W},\mathbf{Z}$) or Greek ($\mathbf{\Lambda}$) letters. We specify $\mathbf{J}$ as the all-ones matrix, in which every element is equal to one. The symbols $\odot$ and $\oslash$ represent the Hadamard product and Hadamard division, respectively. $\|\cdot\|_{\mathrm{F}}$ is the Frobenius norm of a matrix. $\langle \cdot,\cdot \rangle_{\mathrm{F}}$ is Frobenius inner product of two matrices.

\subsection{Conventional Adaptive Smoothing Method (ASM)}

Denote by $\mathbf{Z}(x,t)$ the macroscopic traffic speed field with space and time indices $x$ and $t$. Note that $\mathbf{Z}(x,t)$ can also analogously denote the traffic flux or density field. Given a set of traffic speed measurements $\{(x_n,t_n,v_n)\}_{n}$, ASM first calculates two (a priori) speed fields,
\begin{equation} 
    \label{eqn:fields}
    \begin{aligned}
     &\mathbf{Z}^{\rm free} (x,t) = \frac{1}{N(x,t)} \sum_n \phi \left( x-x_n, t-t_n-\frac{x-x_n}{c_{\rm free}} \right) v_n, \\
     &\mathbf{Z}^{\rm cong} (x,t) = \frac{1}{N(x,t)} \sum_n \phi \left( x-x_n, t-t_n-\frac{x-x_n}{c_{\rm cong}} \right) v_n,
    \end{aligned}
\end{equation}
corresponding to two wave speeds $c_{\rm free}$ and $c_{\rm cong}$, where $N(x,t)$ is a normalization constant and $\phi(\cdot,\cdot)$ is a kernel function (e.g., bi-variate Gaussian).

The ASM defines $\mathbf{Z}(x,t)$ using the following convex combination, written in a compact form as
\begin{equation}
    \begin{aligned}
     \mathbf{Z} = \mathbf{W}^{\rm cong} \odot \mathbf{Z}^{\rm cong} + (\mathbf{J}-\mathbf{W}^{\rm cong}) \odot \mathbf{Z}^{\rm free}
    \end{aligned}
\end{equation}
with the weight field $\mathbf{W}^{\rm cong}(x,t) \in [0,1]$ defined as:
\begin{equation}
\label{eqn:pre-weight}
\begin{aligned}
\mathbf{W}^{\rm cong} = \frac{1}{2} \left[ 1+\tanh\left( \frac{V_{\rm thr}-\min \{ \mathbf{Z}^{\rm free}, \mathbf{Z}^{\rm cong} \} }{\Delta V} \right) \right],
\end{aligned}
\end{equation}
where $V_{\rm thr}$ is the threshold speed and $\Delta V$ is the transition width, which depends on field traffic conditions and requires independent tuning. Note the operators $\min\{\cdot,\cdot\}$ and $\tanh(\cdot)$ in \eqref{eqn:pre-weight} are applied element-wise.

The ASM weight field \eqref{eqn:pre-weight} is based on the observation that propagating structures in congested traffic are very persistent, so that it favours the congested estimation if any of the two a priori estimates indicates congested traffic. It is a heuristic formula that might not adapt well to free-flow and congested traffic speeds, and is sensitive at near-capacity traffic conditions. Also, the dependence on field parameters $V_{\rm thr}$ and $\Delta V$ requires frequent calibration to handle dynamic traffic conditions. Further, $\eqref{eqn:pre-weight}$ is only defined for two traffic wave speeds, which is inadequate when reproducing heterogeneous traffic dynamics (wave speeds can be any convex combination of these two extreme speeds). 

To this end, we formulate the weight calculation $\eqref{eqn:pre-weight}$ as a constrained optimization problem that adapts to traffic conditions in the observed data and eliminates the dependence on field parameters. Our framework can also accommodate multiple wave speeds, which allows one to reproduce richer traffic dynamics. This is shown next.

\subsection{Optimal ASM Weight Calculation}

We formulate the weight calculation as the following matrix completion problem:
\begin{equation} \label{eqn:opt1}
\begin{aligned}
    \underset{\{\mathbf{W}^i\}_{i=1}^m}{\mathrm{Minimize}} \quad & \Big\|\mathsf{P}_\Omega \Big( \sum_i (\mathbf{W}^i \odot\mathbf{Z}^i) - \mathbf{Z}\Big) \Big\|_{\mathrm{F}} \\
    \textrm{s.t.} \quad & \sum_i \mathbf{W}^i = \mathbf{J} ,
\end{aligned}
\end{equation}
where $\mathbf{Z}^i$ is the $i$th a priori speed field estimate corresponding to a wave speed $c_i$, determined using \eqref{eqn:fields}. The associated weights $\mathbf{W}^i$ are the decision variables of \eqref{eqn:opt1}. The binary mask operator $\mathsf{P}_\Omega$ evaluates the objective function only at the observed indices $\Omega$.

For a given set of wave speeds $c_i$ and priori estimate $\mathbf{Z}^i$, the problem \eqref{eqn:opt1} determines the optimal convex coefficients $\mathbf{W}^i$ that results in least norm error. For $m=2$, \eqref{eqn:opt1} solves for weights of two wave speeds, say $c_{\rm free}$ and $c_{\rm cong}$, similar to the ASM setting. Naturally, \eqref{eqn:opt1} can easily consider multiple wave speeds $m>2$ in the final speed field estimation instead of just two wave speeds as in the ASM. 

\subsection{Solving \eqref{eqn:opt1} using Alternating Direction Method of Multipliers}

We apply Alternating Direction Method of Multipliers (ADMM) \cite{boyd2011distributed} to solve \eqref{eqn:opt1}. Using an auxiliary variable $\hat{\mathbf{Z}}$, we first reformulate $\eqref{eqn:opt1}$ as:
\begin{equation} \label{eqn:opt2}
\begin{aligned}
    \underset{\{\mathbf{W}^i\}_{i=1}^m,\hat{\mathbf{Z}}}{\mathrm{Minimize}} \quad & \|\mathsf{P}_\Omega(\hat{\mathbf{Z}} - \mathbf{Z})\|_{\mathrm{F}} \\
    \textrm{s.t.} \quad & \hat{\mathbf{Z}} = \sum_i \mathbf{W}^i \odot\mathbf{Z}^i  \\
    \quad \quad & \sum_i \mathbf{W}^i = \mathbf{J}
\end{aligned}
\end{equation}

With the introduction of $\hat{\mathbf{Z}}$, we can simplify the objective function that was hard to directly minimize. Based on this reformation, we write the augmented Lagrangian function of \eqref{eqn:opt2} as follows:
\begin{equation}
\begin{aligned}
& L_\beta(\hat{\mathbf{Z}},\mathbf{W}^1,...,\mathbf{W}^m,\mathbf{\Lambda}_1,\mathbf{\Lambda}_2) = \frac{1}{2}\big\|\mathbf{M}\odot(\mathbf{Z}-\hat{\mathbf{Z}})\big\|^2_{\mathrm{F}} \\
&\quad +\big\langle\mathbf{\Lambda}_1,\hat{\mathbf{Z}}-\sum_i \mathbf{W}^i\odot \mathbf{Z}^i\big\rangle_{\mathrm{F}}+\frac{\beta}{2}\big\|\hat{\mathbf{Z}}-\sum_i \mathbf{W}^i\odot \mathbf{Z}^i\big\|_{\mathrm{F}}^2 \\
&\quad +\big\langle\mathbf{\Lambda}_2,\sum_i \mathbf{W}^i-\mathbf{J}\big\rangle_{\mathrm{F}}+\frac{\beta}{2}\big\|\sum_i \mathbf{W}^i-\mathbf{J}\big\|_{\mathrm{F}}^2,
\end{aligned}
\end{equation}
where $\hat{\mathbf{Z}},\mathbf{W}^1,...,\mathbf{W}^m$ are the primal variables and $\mathbf{\Lambda}_1$, $\mathbf{\Lambda}_2$ are the dual variables for \eqref{eqn:opt2}.

The ADMM then proceeds iteratively by alternatively optimizing the primal variables and dual variables. This is summarized in Algorithm~\ref{algo1}.


{ 
\small
\begin{algorithm}[h!]
\label{algo1}
\caption{ADMM for solving problem \eqref{eqn:opt2}}
\KwData{$\mathbf{Z}$}
\KwResult{$\mathbf{W}^1,...,\mathbf{W}^m$}

\textbf{Initialize}: $\hat{\mathbf{Z}},\mathbf{W}^1,...,\mathbf{W}^m,\mathbf{\Lambda}_1,\mathbf{\Lambda}_2$

\While{stopping criterion not met}{
    \begin{align*}
    \hat{\mathbf{Z}} &\mapsfrom \underset{\hat{\mathbf{Z}}}{\arg\min} \Big(\frac{1}{2}\big\|\mathbf{M}\odot(\mathbf{Z}-\hat{\mathbf{Z}})\big\|^2_{\mathrm{F}} \\
    &+ \big\langle\mathbf{\Lambda}_1,\hat{\mathbf{Z}}-\sum_i \mathbf{W}^i\odot \mathbf{Z}^i\big\rangle_{\mathrm{F}} \\
    &+ \frac{\beta}{2}\big\|\hat{\mathbf{Z}}-\sum_i \mathbf{W}^i\odot \mathbf{Z}^i\big\|_{\mathrm{F}}^2\Big) \\
    \mathbf{W}^i &\mapsfrom \underset{\mathbf{W}^i}{\arg\min} \Big(\big\langle\mathbf{\Lambda}_1,\hat{\mathbf{Z}}-\sum_i \mathbf{W}^i\odot \mathbf{Z}^i\big\rangle_{\mathrm{F}} \\
    &+ \frac{\beta}{2}\big\|\hat{\mathbf{Z}}-\sum_i \mathbf{W}^i\odot \mathbf{Z}^i\big\|_{\mathrm{F}}^2 \\
    &+ \big\langle\mathbf{\Lambda}_2,\sum_i \mathbf{W}^i-\mathbf{J}\big\rangle_{\mathrm{F}} + \frac{\beta}{2}\big\|\sum_i \mathbf{W}^i-\mathbf{J}\big\|_{\mathrm{F}}^2\Big)
    \end{align*}

    
    \quad $\mathbf{\Lambda}_1 \mapsfrom \mathbf{\Lambda}_1 + \beta(\hat{\mathbf{Z}}-\sum_i \mathbf{W}^i\odot \mathbf{Z}^i)$
    
    \quad $\mathbf{\Lambda}_2 \mapsfrom \mathbf{\Lambda}_2 + \beta(\sum_i \mathbf{W}^i-\mathbf{J})$
}
\end{algorithm}
}

The primal variable updates in Algorithm \ref{algo1} can be solved in closed form:
\begin{equation} \nonumber
    \begin{aligned}
     &\hat{\mathbf{Z}} = (\mathbf{M}\odot \mathbf{Z}-\mathbf{\Lambda}_1+\beta(\sum_i \mathbf{W}^i\odot \mathbf{Z}^i)) \oslash (\mathbf{M}+\beta \mathbf{J}) \\
     &\mathbf{W}^i = (\beta(\hat{\mathbf{Z}}\odot \mathbf{Z}^i+\mathbf{J}-\sum_{r \neq i}\mathbf{W}^r\odot(\mathbf{Z}^r\odot \mathbf{Z}^i+\mathbf{J})) \\
        & \quad \quad \quad + \mathbf{\Lambda}_1\odot \mathbf{Z}^i -\mathbf{\Lambda}_1) 
            \oslash (\beta(\mathbf{Z}^i\odot \mathbf{Z}^i+\mathbf{J}))
    \end{aligned}
\end{equation}
The dual variables are updated using simple gradient ascent. Algorithm~\ref{algo1} terminates when the primal and dual residuals reach their pre-defined feasibility tolerances \cite{boyd2011distributed}.







\section{EXPERIMENTAL SETUP}
\label{sec4:experiments}

\begin{figure}[thpb]
  \centering
  \includegraphics[scale=0.75]{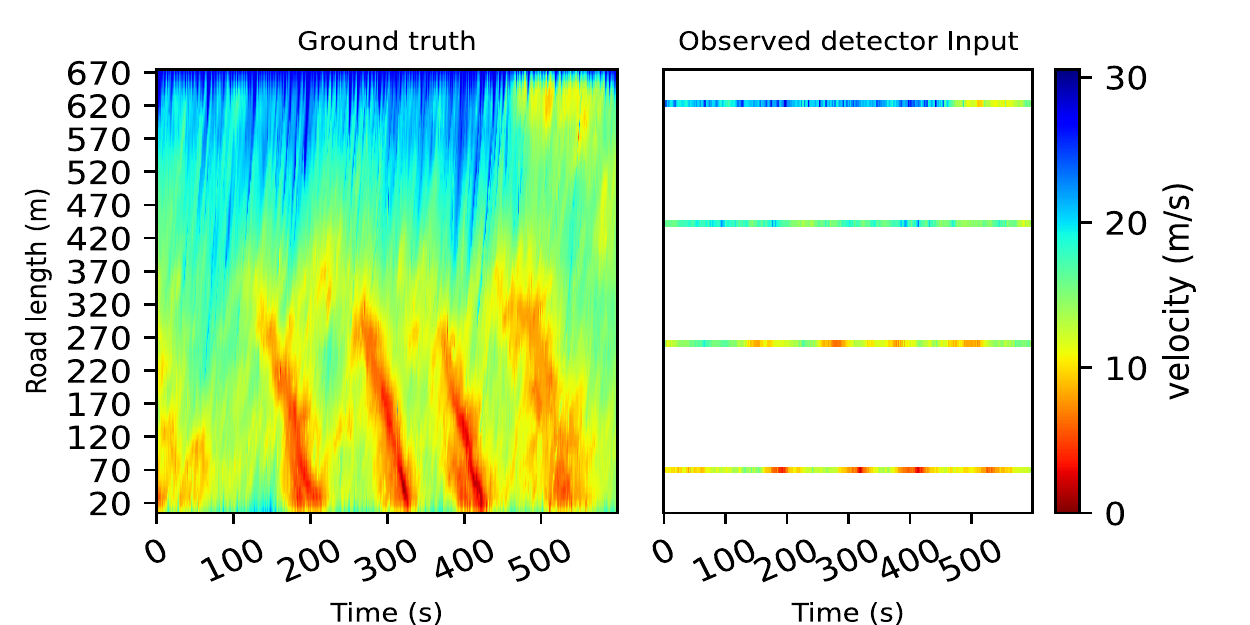}
  (a) \hspace{1.0 in} (b)
  \caption{Case 1 data (NGSIM US 101 lane 2, $100-700$ s): (a) Ground truth data and (b) Input measurements (4 detectors).}
  \label{figure1}
\end{figure}

\begin{figure}[thpb]
  \centering
  \includegraphics[scale=0.75]{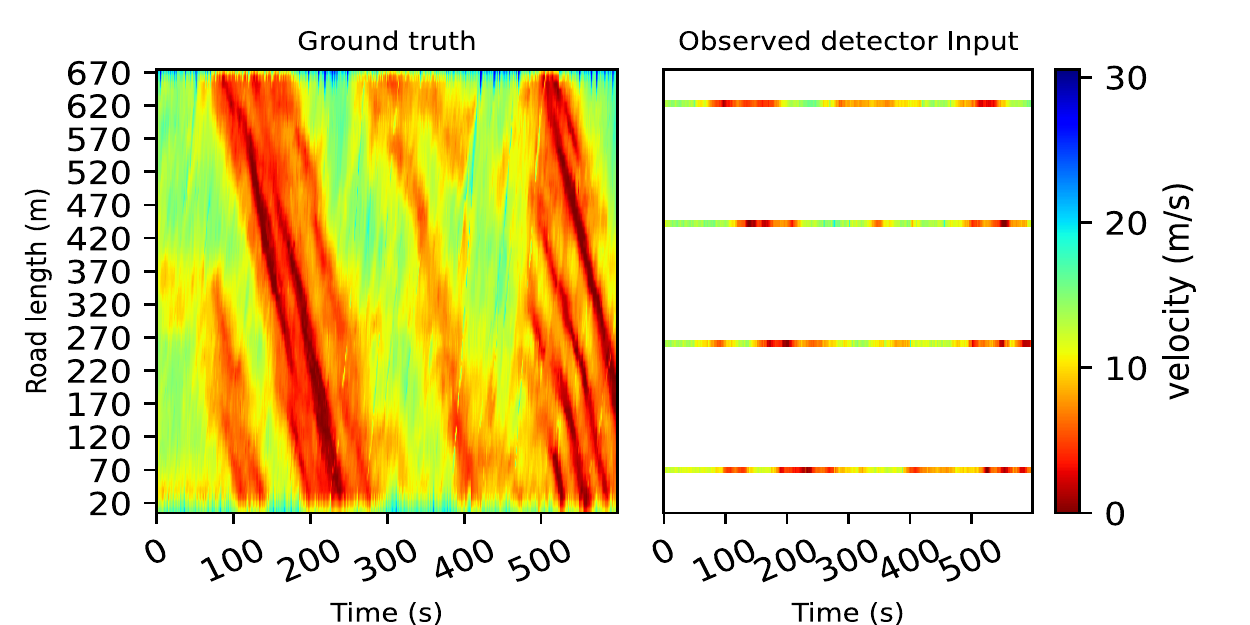}
  (a) \hspace{1.0 in} (b)
  \caption{Case 2 data (NGSIM US 101 lane 2, $1400-2000$ s): (a) Ground truth data and (b) Input measurements (4 detectors).}
  \label{figure2}
\end{figure}

We test the estimation performance of our proposed algorithm using the NGSIM US 101 highway data \cite{ngsim}. The data consists of vehicle trajectories from a road section $670$ m in length over a $2700$ s time period. This data is converted to a ground-truth macroscopic speed field $\mathbf{Z}$ of dimension $67 \times 2700$. Heatmaps of these data are shown in Fig.~\ref{figure1}a and Fig.~\ref{figure2}a for $100-700$ s and and $1400-2200$ s as two cases. The $100-700$ s case is a mixture of free flow and congestion while $1400-2200$ s consists mostly of congested traffic (backward shockwaves). We assume that the input measurements come from  stationary detector data, though the methodology can also be applied to floating car data. Thus, the set of input indices $\Omega$ consists of the number of rows of $\mathbf{Z}$ observed.

We conduct three experiments in this study. In the first experiment, we determine optimal weights for two wave speeds (i.e., $m=2$), and compare the estimation error of the proposed ADMM-based algorithm and the conventional ASM. We use the wave speeds $c_1=80$ km/hr and $c_2=-15$ km/hr, corresponding to the optimal setting for ASM. The estimation methods are compared for two different time periods to evaluate their performance in free-flowing and congested traffic conditions. The true speed fields and input measurements for these two cases are shown in Fig.~\ref{figure1}b and Fig.~\ref{figure2}b.

In the second experiment, we evaluate the benefits of the proposed ADMM-based estimation algorithm to incorporate multiple wave speeds (i.e., $m>2$). We consider the following set of wave speeds: $\{ -20, -17.5, -15, -12.5, +60, +70, +80, +90 \}$ km/hr. These wave speeds are chosen so as to incorporate a range of traffic waves emanating in free-flowing and congested traffic conditions. In the third experiment, we investigate our method's performance on different level of input sparsity by changing the detector coverage rate.


The parameters of ASM and ADMM used in our experiments are summarized in Table \ref{table_parameter}. 
The spatial and temporal smoothing widths, $\sigma$ and  $\tau$, are the parameters of the kernel function $\phi(\cdot,\cdot)$. We note that the choice of $\sigma$ and $\tau$ are sensitive to estimation errors, especially for lower reconstruction window. We chose them as half the average inter-detector spacing and sampling time, respectively, as recommended in \cite{trieber2011filter}. $\beta$ is an ADMM hyper-parameter that controls the convergence rate; smaller values of $\beta$ means larger step sizes and faster convergence, but this can cause instability \cite{sun2014alternating}. We use $\beta=1$ in our experiments.

\begin{table}[h!]
\small
\caption{Parameter Setting}
\label{table_parameter}
\begin{center}
\begin{tabular}{|c|c|l|}
\hline
Parameter & Value & Description\\
\hline
$v_{\rm thr}$ & 60 km/h & Critical traffic speed\\
$\Delta v$ & 20 km/h & Transition width\\
$\sigma$ & $\Delta x/2$ & Space coordinate smoothing width\\
$\tau$ & $\Delta t/2$ & Time coordinate smoothing width\\
$\beta$ & 1 & Step-size in ADMM\\
\hline
\end{tabular}
\end{center}
\end{table}

The quality of the estimated speed field $\hat{\mathbf{Z}}$ is measured using the relative error in all the experiments,
\begin{equation*}
    m_{\mathrm{r}} = \frac{\|\hat{\mathbf{Z}} - \mathbf{Z}\|_{\mathrm{F}}}{\|\mathbf{Z}\|_{\mathrm{F}}}
\end{equation*}



\section{Result Analysis and Discussion}
\label{sec5:results}

The experiment results are discussed below: 
\newline

\paragraph{Expt 1: Comparison of ASM and ADMM-based} The estimation errors $m_r$ for different cases shown in Fig.~\ref{figure1} and Fig.~\ref{figure2} are summarized in Table \ref{table_result}. The least error is highlighted in bold. We see that the ADMM estimation gives slightly better performance in both cases, which speaks to the utility of the weights estimation using ADMM (the main difference between ADMM and ASM in this case).
\begin{table}[h!]
\small
\caption{Estimation errors for Expt 1}
\label{table_result}
\begin{center}
\begin{tabular}{|c|c||c|c|}
\hline
Case 1 & $m_{\mathrm{r}}$  & Case 2 & $m_{\mathrm{r}}$ \\
\hline
ASM & 0.12417 & ASM & 0.19843 \\
ADMM & \textbf{0.12054} & ADMM & \textbf{0.19637} \\
\hline
\end{tabular}
\end{center}
\end{table}

\begin{figure}[h!]
      \centering
      \includegraphics[scale=0.75]{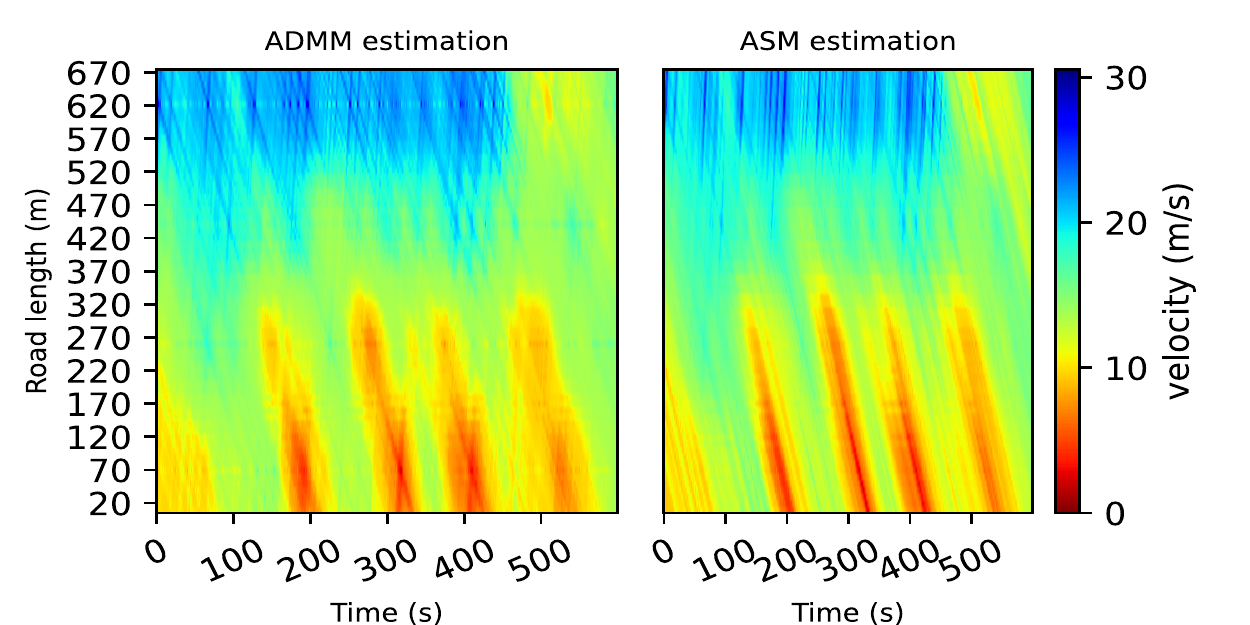}
      \caption{Comparison of the ADMM-based and ASM estimation results for the Case 1 data (Refer Fig.~\ref{figure1}).}
      \label{figureRC1}
   \end{figure}

\begin{figure}[h!]
      \centering
      \includegraphics[scale=0.75]{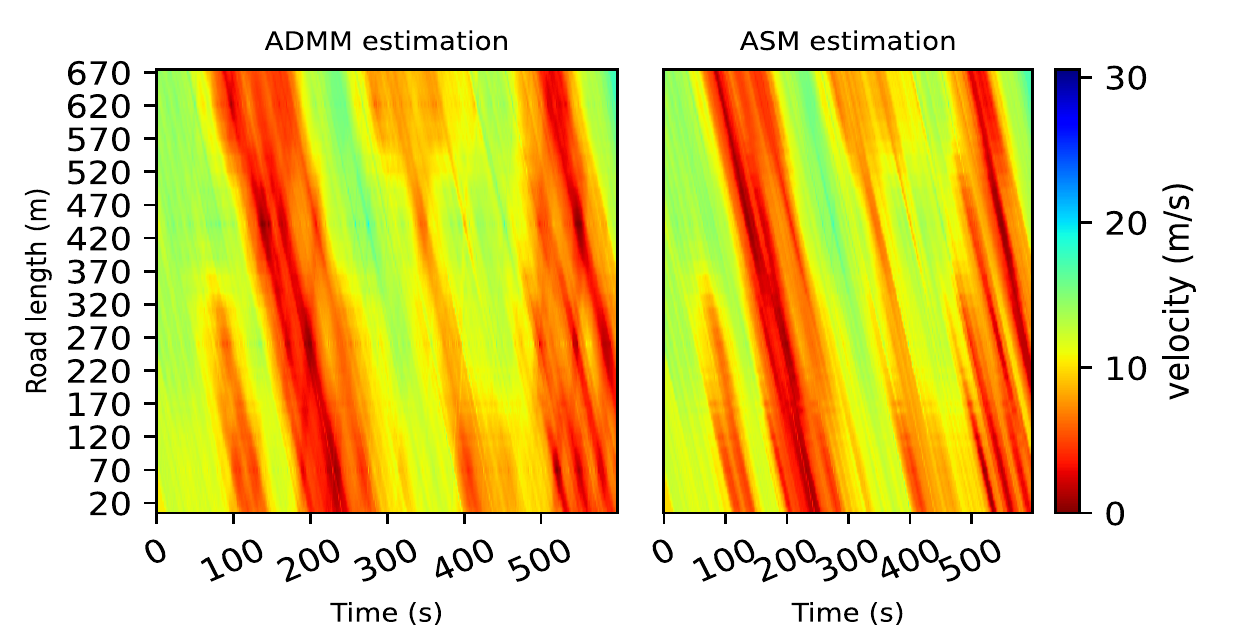}
      \caption{Comparison of the ADMM-based and ASM estimation results for the Case 2 data (Refer Fig.~\ref{figure2}).}
      \label{figureRC2}
   \end{figure}

We also show the visualization for both method's estimation result in Fig. \ref{figureRC1} and Fig. \ref{figureRC2}. Both algorithms successfully captured the ground truth wave propagation dynamics. It is also notable that ASM  tends to show a thinner and sharper pattern than ADMM, which is not an accurate representation for ground truth. This is more obvious in the estimation shown in Fig. \ref{figureRC2}. In free flow region, we see that ADMM's estimation consist of several minor waves, which is not observed in the ASM's estimation. Nevertheless, ADMM captured the area size and shape more accurately than ASM.

\begin{figure}[h!]
      \centering
      \includegraphics[scale=0.65]{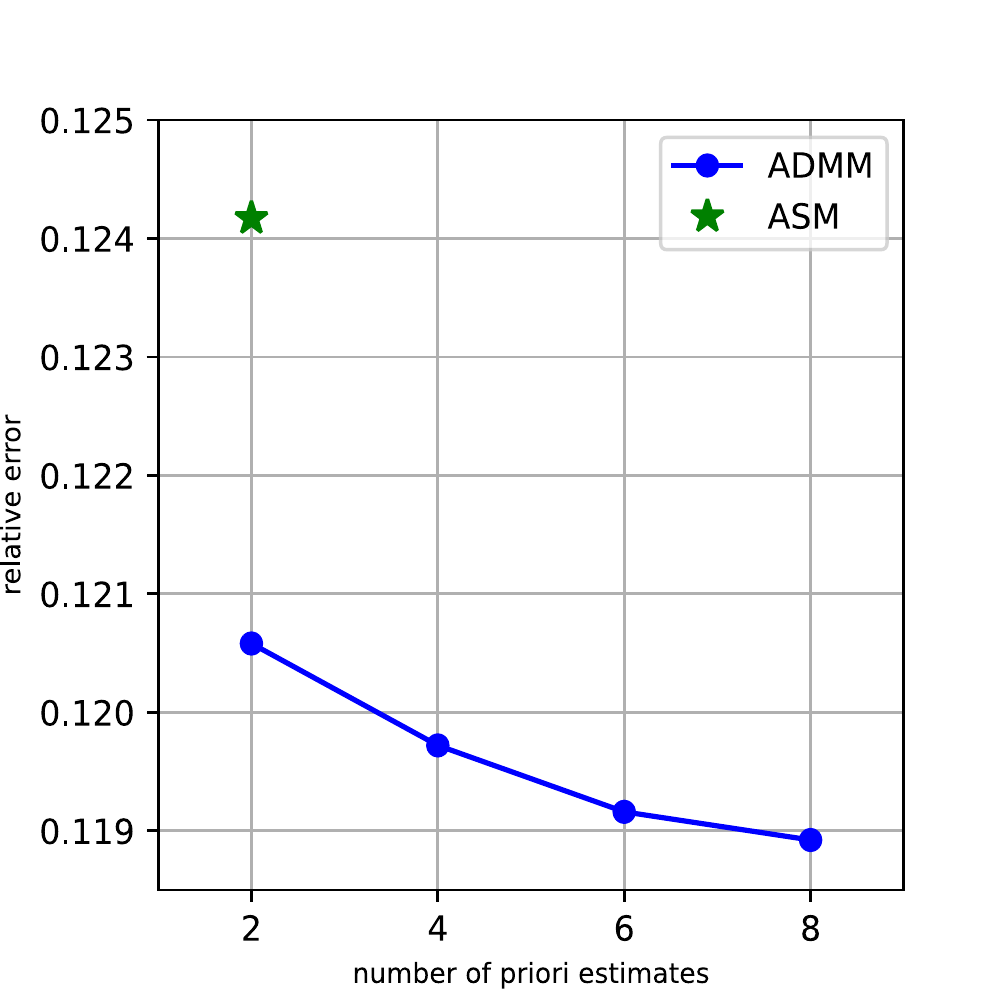}
      \caption{Error rates for the ADMM-based algorithm considering multiple a priori estimates (i.e., with multiple wave speeds). The results corresponds to the input data shown in Fig.~\ref{figure2}.}
      \label{figureM}
   \end{figure}

\paragraph{Expt 2: Multiple a priori estimates}
Next, we evaluate the ADMM estimation error when considering multiple wave speeds (a priori estimates). We consider the wave speeds: $[-20,90;-17.5,80;-15,70;-12.5,65]$. To measure the improvement, we run different sub-experiments, where in each experiment, we add one pair of congested and free-flow prior estimates and calculate the our algorithm's estimation error. For example, we start with two a priori estimates with $c_{\rm cong} = -10$ and $c_{\rm free} = 90$, then add next pair of wave speeds $c_{\rm cong} = -17.5$ and $c_{\rm free} = 80$, and so on. The results are summarized in Fig. \ref{figureM}. The green star marker represents the estimation error of conventional ASM with two a priori estimates for reference. We observe that increasing the number of wave speeds could further reduce reconstruction error as seen in Fig. \ref{figureM}. The choice of wave speeds considered in this experiment is rather heuristic. A more careful choice of wave speeds could lead to even greater improvement in the performance.

\begin{figure}[h!]
  \centering
  \includegraphics[scale=0.5]{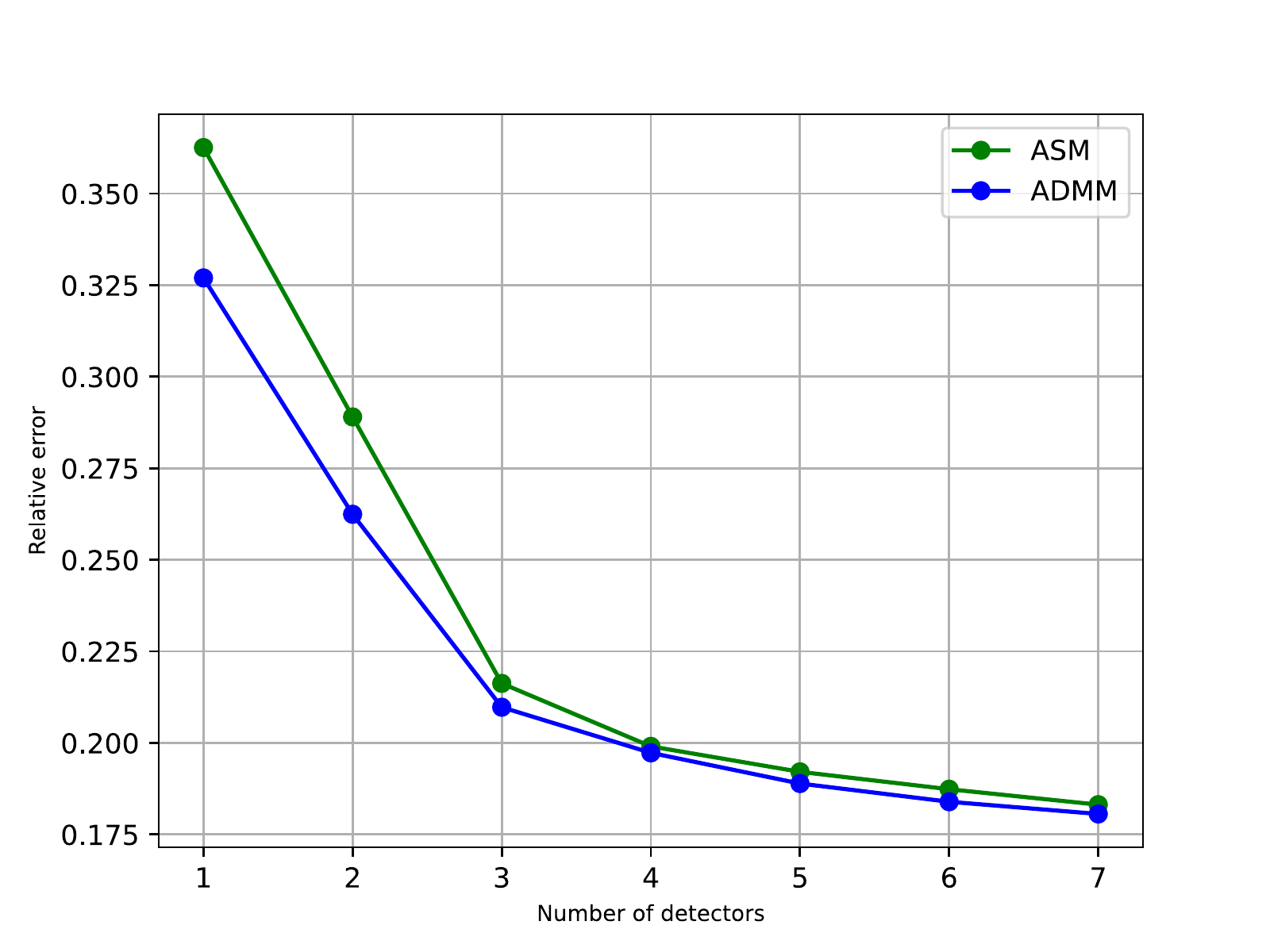}
  \caption{Error rates comparison for the ADMM-based and ASM algorithm at different coverage rates (i.e., number of detectors).}
  \label{figureCR}
\end{figure}

\paragraph{Expt 3: Different input coverage rates}
We test our proposed algorithm's performance under different detector coverage rates. As shown in Fig. \ref{figureCR}, we applied the ADMM-based (with $m=2$) and ASM algorithms on the case 2 data for different number of detectors (i.e., input penetration rates). The curves in Fig.~\ref{figureCR} show the error rate trend for increasing the input coverage rates from $1\%$ to $10\%$ (i.e., number of detectors from $1$ to $7$). We observe that the errors from both the algorithms decrease with the increase of input information, with ADMM performing better than ASM.
But, the difference in error rates are higher when the input information is very limited, i.e., $1\%-5\%$. This indicates that our ADMM-based algorithm extracts more information from the sparse input measurements in comparison to the conventional ASM.
   
\section{CONCLUSIONS}
\label{sec6:conclusions}

This study improves the conventional
Adaptive Smoothing Method (ASM) for traffic state estimation by proposing a systematic procedure to calculate the weights of congested and free-flow traffic speed fields (a priori estimates), which is otherwise done heuristically. We formulate the a priori estimates' weight calculation as a matrix completion problem, and efficiently solve it using the Alternating Direction Method of Multipliers algorithm. Our algorithm doesn't depend on field parameters, and can, thus, reduce the effort involved in field calibrations as required in conventional ASM. Our framework allows one to consider a priori speed field estimates corresponding to multiple wave speeds rather than just two wave speeds as done in the conventional ASM. This is advantageous for reproducing richer traffic dynamics, e.g., vehicular traffic flows with wide range of desired speed distributions. 
Experiments using real traffic data show that our proposed algorithm reduces the estimation error and achieves better performance than the conventional ASM, particularly when using multiple a priori input speed field estimates. The experiments also show reduction in estimation error when the input measurements are sparse.

We observe in our multiple a priori estimates experiments that proper choice of wave speed $c_i$ can enhance the performance of our algorithm's performance to a great extent. However, the choice of these values can be rather heuristic. We believe that the present approach can be extended to incorporate tuning these parameters in a more effective manner and improve the estimation accuracy further. We also observe that our algorithm tends to estimate high velocity regions more accurately, which is a natural outcome due to the objective function based on Frobenius error minimization, which tends to focus on high magnitude numbers. In the context of traffic estimation where the congested region is of greater interest, one can utilize a different objective function, which emphasizes errors in low-speed traffic or by changing the estimated traffic state variable from velocity to density. We leave this to future research.

\section*{Acknowledgment}
This work was supported by the NYUAD Center for Interacting Urban Networks (CITIES), funded by Tamkeen under the NYUAD Research Institute Award CG001.  The opinions expressed in this article are those of the authors alone do not represent the opinions of CITIES.

	
	
{ \small
\bibliographystyle{plainnat}
\bibliography{References}

\begin{thebibliography}{29}
\providecommand{\natexlab}[1]{#1}
\providecommand{\url}[1]{\texttt{#1}}
\expandafter\ifx\csname urlstyle\endcsname\relax
  \providecommand{\doi}[1]{doi: #1}\else
  \providecommand{\doi}{doi: \begingroup \urlstyle{rm}\Url}\fi

\bibitem[{Bekiaris-Liberis} et~al.(2016){Bekiaris-Liberis}, {Roncoli}, and
  {Papageorgiou}]{bekiaris2016tse_ieee}
N.~{Bekiaris-Liberis}, C.~{Roncoli}, and M.~{Papageorgiou}.
\newblock Highway traffic state estimation with mixed connected and
  conventional vehicles.
\newblock \emph{IEEE Transactions on Intelligent Transportation Systems},
  17\penalty0 (12):\penalty0 3484--3497, 2016.
\newblock \doi{10.1109/TITS.2016.2552639}.

\bibitem[Benkraouda et~al.(2020)Benkraouda, Thodi, Yeo, Menendez, and
  Jabari]{benkraouda2020traffic}
Ouafa Benkraouda, Bilal~Thonnam Thodi, Hwasoo Yeo, Monica Menendez, and
  Saif~Eddin Jabari.
\newblock Traffic data imputation using deep convolutional neural networks.
\newblock \emph{IEEE Access}, 8:\penalty0 104740--104752, 2020.

\bibitem[Boyd et~al.(2011)Boyd, Parikh, Chu, Peleato, Eckstein,
  et~al.]{boyd2011distributed}
Stephen Boyd, Neal Parikh, Eric Chu, Borja Peleato, Jonathan Eckstein, et~al.
\newblock Distributed optimization and statistical learning via the alternating
  direction method of multipliers.
\newblock \emph{Foundations and Trends{\textregistered} in Machine learning},
  3\penalty0 (1):\penalty0 1--122, 2011.

\bibitem[Chen et~al.(2019)Chen, Zhang, Li, and Li]{chen2019filter}
X.~Chen, S.~Zhang, L.~Li, and L.~Li.
\newblock Adaptive rolling smoothing with heterogeneous data for traffic state
  estimation and prediction.
\newblock \emph{IEEE Transactions on Intelligent Transportation Systems},
  20\penalty0 (4):\penalty0 1247--1258, 2019.

\bibitem[Huang and Agarwal(2020)]{huang2020physics}
Jiheng Huang and Shaurya Agarwal.
\newblock Physics informed deep learning for traffic state estimation.
\newblock In \emph{IEEE Intelligent Transportation Systems Conference}, 2020.

\bibitem[Jabari and Liu(2012)]{jabari2012stochastic}
Saif~Eddin Jabari and H.~Liu.
\newblock A stochastic model of traffic flow: {T}heoretical foundations.
\newblock \emph{Transportation Research Part B: Methodological}, 46\penalty0
  (1):\penalty0 156--174, 2012.

\bibitem[Jabari and Liu(2013)]{jabari2013gauss}
Saif~Eddin Jabari and H.~Liu.
\newblock A stochastic model of traffic flow: {G}aussian approximation and
  estimation.
\newblock \emph{Transportation Research Part B: Methodological}, 47:\penalty0
  15--41, 2013.

\bibitem[Jabari et~al.(2018)Jabari, Zheng, Liu, and
  Filipovska]{jabari2018stochastic}
Saif~Eddin Jabari, F.~Zheng, H.~Liu, and M.~Filipovska.
\newblock Stochastic {L}agrangian modeling of traffic dynamics.
\newblock In \emph{The 97th Annual Meeting of the Transportation Research
  Board, Washington D.C}, pages 18--04170, 2018.

\bibitem[Jabari et~al.(2019)Jabari, Dilip, Lin, and
  Thonnam~Thodi]{jabari2019learning}
Saif~Eddin Jabari, D.~Dilip, D.~Lin, and B.~Thonnam~Thodi.
\newblock Learning traffic flow dynamics using random fields.
\newblock \emph{IEEE Access}, 7:\penalty0 130566--130577, 2019.

\bibitem[Jabari et~al.(2020)Jabari, Freris, and Dilip]{jabarisparse2020}
Saif~Eddin Jabari, Nikolaos~M. Freris, and Deepthi~Mary Dilip.
\newblock Sparse travel time estimation from streaming data.
\newblock \emph{Transportation Science}, 54\penalty0 (1):\penalty0 1--20, 2020.
\newblock \doi{10.1287/trsc.2019.0920}.
\newblock URL \url{https://doi.org/10.1287/trsc.2019.0920}.

\bibitem[Jia et~al.(2016)Jia, Wu, and Du]{yuhan2016dnn-speed}
Y.~Jia, J.~Wu, and Y.~Du.
\newblock Traffic speed prediction using deep learning method.
\newblock In \emph{2016 IEEE 19th International Conference on Intelligent
  Transportation Systems (ITSC)}, pages 1217--1222, 2016.

\bibitem[Li and Jabari(2019)]{li2019position}
L.~Li and S.E. Jabari.
\newblock Position weighted backpressure intersection control for urban
  networks.
\newblock \emph{Transportation Research Part B: Methodological}, 128:\penalty0
  435--461, 2019.

\bibitem[Li et~al.(2021)Li, Okoth, and Jabari]{li2021backpressure}
Li~Li, Victor Okoth, and Saif~Eddin Jabari.
\newblock Backpressure control with estimated queue lengths for urban network
  traffic.
\newblock \emph{IET Intelligent Transport Systems}, in press:\penalty0 DOI:
  10.1049/itr2.12027, 2021.

\bibitem[Li et~al.(2020)Li, Yang, and Jabari]{li2020short}
Wenqing Li, Chuhan Yang, and Saif~Eddin Jabari.
\newblock Short-term traffic forecasting using high-resolution traffic data.
\newblock In \emph{2020 IEEE 23rd International Conference on Intelligent
  Transportation Systems (ITSC)}, pages 1--6. IEEE, 2020.

\bibitem[Li et~al.(2022)Li, Yang, and Jabari]{li2022ensemblemc}
Wenqing Li, Chuhan Yang, and Saif~Eddin Jabari.
\newblock Nonlinear traffic prediction as a matrix completion problem with
  ensemble learning.
\newblock \emph{Transportation Science}, 56\penalty0 (1):\penalty0 52--78,
  2022.
\newblock \doi{10.1287/trsc.2021.1086}.
\newblock URL \url{https://doi.org/10.1287/trsc.2021.1086}.

\bibitem[Newell(1993)]{newell1993simplified}
Gordon~F Newell.
\newblock A simplified theory of kinematic waves in highway traffic, part i:
  General theory.
\newblock \emph{Transportation Research Part B: Methodological}, 27\penalty0
  (4):\penalty0 281--287, 1993.

\bibitem[Papadopoulou et~al.(2018)Papadopoulou, Roncoli, Bekiaris-Liberis,
  Papamichail, and Papageorgiou]{papadopoulou2018microscopic}
Sofia Papadopoulou, Claudio Roncoli, Nikolaos Bekiaris-Liberis, Ioannis
  Papamichail, and Markos Papageorgiou.
\newblock Microscopic simulation-based validation of a per-lane traffic state
  estimation scheme for highways with connected vehicles.
\newblock \emph{Transportation Research Part C: Emerging Technologies},
  86:\penalty0 441--452, 2018.

\bibitem[Schreiter et~al.(2010)Schreiter, van Lint, Treiber, and
  Hoogendoorn]{schreiter2010two}
Thomas Schreiter, Hans van Lint, Martin Treiber, and Serge Hoogendoorn.
\newblock Two fast implementations of the adaptive smoothing method used in
  highway traffic state estimation.
\newblock In \emph{13th International IEEE Conference on Intelligent
  Transportation Systems}, pages 1202--1208. IEEE, 2010.

\bibitem[Shi et~al.(2021)Shi, Mo, Huang, Di, and Du]{shi2021physics}
Rongye Shi, Zhaobin Mo, Kuang Huang, Xuan Di, and Qiang Du.
\newblock A physics-informed deep learning paradigm for traffic state and
  fundamental diagram estimation.
\newblock \emph{IEEE Transactions on Intelligent Transportation Systems}, 2021.

\bibitem[Sun and Fevotte(2014)]{sun2014alternating}
Dennis~L Sun and Cedric Fevotte.
\newblock Alternating direction method of multipliers for non-negative matrix
  factorization with the beta-divergence.
\newblock In \emph{2014 IEEE international conference on acoustics, speech and
  signal processing (ICASSP)}, pages 6201--6205. IEEE, 2014.

\bibitem[Thodi et~al.(2021)Thodi, Khan, Jabari, and Menéndez]{thodi2021itsc}
Bilal~Thonnam Thodi, Zaid~Saeed Khan, Saif~Eddin Jabari, and Mónica Menéndez.
\newblock Learning traffic speed dynamics from visualizations.
\newblock In \emph{2021 IEEE International Intelligent Transportation Systems
  Conference (ITSC)}, pages 1239--1244, 2021.
\newblock \doi{10.1109/ITSC48978.2021.9564541}.

\bibitem[Thodi et~al.(2022)Thodi, Khan, Jabari, and
  Menendez]{thodi2021anisocnn}
Bilal~Thonnam Thodi, Zaid~Saeed Khan, Saif~Eddin Jabari, and Monica Menendez.
\newblock Incorporating kinematic wave theory into a deep learning method for
  high-resolution traffic speed estimation.
\newblock \emph{IEEE Transactions on Intelligent Transportation Systems}, In
  press, 2022.

\bibitem[Treiber and Helbing(2002)]{treiber2002filter}
M.~Treiber and D.~Helbing.
\newblock Reconstructing the spatio-temporal traffic dynamics from stationary
  detector data.
\newblock \emph{Cooperative Transportation Dynamics}, 1\penalty0 (3):\penalty0
  3--1, 2002.

\bibitem[Treiber et~al.(2011)Treiber, Kesting, and Wilson]{trieber2011filter}
Martin Treiber, Arne Kesting, and R.~Eddie Wilson.
\newblock Reconstructing the traffic state by fusion of heterogeneous data.
\newblock \emph{Computer-Aided Civil and Infrastructure Engineering},
  26\penalty0 (6):\penalty0 408--419, 2011.

\bibitem[{United States Department of Transportation}(2006)]{ngsim}
{United States Department of Transportation}.
\newblock {NGSIM—Next Generation Simulation}, 2006.
\newblock URL \url{https://ops.fhwa.dot.gov/trafficanalysistools/ngsim.htm}.

\bibitem[Xiao et~al.(2018)Xiao, Wei, and Liu]{xiao2018speed}
Jianli Xiao, Chao Wei, and Yuncai Liu.
\newblock Speed estimation of traffic flow using multiple kernel support vector
  regression.
\newblock \emph{Physica A: Statistical Mechanics and its Applications},
  509:\penalty0 989--997, 2018.

\bibitem[Yang and Menendez(2019)]{kaidi2019queueest}
K.~Yang and M.~Menendez.
\newblock Queue estimation in a connected vehicle environment: {A} convex
  approach.
\newblock \emph{IEEE Transactions on Intelligent Transportation Systems},
  20\penalty0 (7):\penalty0 2480--2496, 2019.

\bibitem[Yang et~al.(2016)Yang, Guler, and Menendez]{yang2016cav}
Kaidi Yang, S.~Ilgin Guler, and Monica Menendez.
\newblock Isolated intersection control for various levels of vehicle
  technology: Conventional, connected, and automated vehicles.
\newblock \emph{Transportation Research Part C: Emerging Technologies},
  72:\penalty0 109 -- 129, 2016.
\newblock ISSN 0968-090X.
\newblock \doi{https://doi.org/10.1016/j.trc.2016.08.009}.
\newblock URL
  \url{http://www.sciencedirect.com/science/article/pii/S0968090X16301437}.

\bibitem[{Yuan} et~al.(2012){Yuan}, {van Lint}, {Wilson}, {van
  Wageningen-Kessels}, and {Hoogendoorn}]{hoogen2012lagrang}
Y.~{Yuan}, J.~W.~C. {van Lint}, R.~E. {Wilson}, F.~{van Wageningen-Kessels},
  and S.~P. {Hoogendoorn}.
\newblock Real-time lagrangian traffic state estimator for freeways.
\newblock \emph{IEEE Transactions on Intelligent Transportation Systems},
  13\penalty0 (1):\penalty0 59--70, 2012.
\newblock \doi{10.1109/TITS.2011.2178837}.

\end{thebibliography}
}
	
	
	
	
	

\end{document}